# On some interesting tasks fron classical theory of numbers

Zurab Aghdgomelashvili

Doctor of Mathematics

Department of Mathematics of Faculty of Informatics and control
Systems of Georgian Technical University, 77, M. Kostava
Street, 0175, Tbilisi, Georgia

E-mail: z.aghdgomelashvili@gtu.ge

**Abstract:** Expressions of type $(p^q - 1)/(p - 1)$ and $(a^2 + ab + b^2)$, where a, b are natural and $p$, $q$ are prime numbers, are studied.

## Introduction

Solving equations in integers is one of the most beautiful parts of mathematics. Over time, many methods have been collected for solving specific Diophantine equations, but general methods for investigating them were developed only in the XX century. In 1900, at the X International Congress of Mathematicians, held in Paris, the famous German mathematician D. Hilbert formulated 23 key mathematical problems in his report. Hilbert's tenth problem was called as "The Task of solving Diophantine equations". Its essence was as follows: is given a Diophantine equation with integer coefficients and several unknowns. Does exist or not an algorithm with a finite number of operations for solving this equation in integers.

We have shown that the entire class of tasks would be easily solved using one simple lemma; in addition, is stated the task to study $\dfrac{p^q - 1}{p - 1}$ type expressions, where $p$ and $q$ are prime numbers, with properties; some properties of the $(a^2 + ab + b^2)$ type expression were studied, and so on.

The following issues are considered in the work:

– Obviously, for prime numbers $p$ and $q$, it is of great interest to determine the quantity of those prime divisors of the $A = \dfrac{p^q - 1}{p - 1}$ number that are less than $p$. To do this, we have considered:

Theorem 1. Let's say that $p$ and $q$ are odd prime numbers and $p = 2q + 1$. Then from the various prime divisors of number $A=\frac{p^q-1}{p-1}$, taken separately, only one is less than *p*. *A* has at least two different prime divisors.

Theorem 2. Let's say that *p* and *q* are odd prime numbers and $p < 2q + 1$. Then all prime divisors of the number $A=\frac{p^q-1}{p-1}$ are greater than *p*;

Theorem 3. Let's say that *q* is an odd prime number and $p \in N \setminus \{1\}$, $p \in ]1;q] \cup [q+2;2q]$, then each of the various prime divisors of the number $A=\frac{p^q-1}{p-1}$, taken separately, is greater than *p*;

Theorem 4. Let's say that *q* is a prime odd number and $p \in \{q + 1; 2q + 1\}$, then from different prime divisors of number $A=\frac{p^q-1}{p-1}$, taken separately, only one of them is less than *p*. *A* has at least two different prime divisors.

Task 1. Let's solve the equation $2^x=\frac{y^z-1}{y-1}$ in natural *x*, *y*, *z* numbers. At the same time *y* must be a prime number.

Task 2. Let's solve the equation $3^x=\frac{y^z-1}{y-1}$ in natural *x*, *y*, *z* numbers. At the same time *y* must be a prime number.

Task 3. Let's solve the equation $p^x=\frac{y^z-1}{y-1}$, where $p \in \{5; 7; 11; 13; ...\}$ prime number *x*, $y \in N$ and *y* – is prime number.

I – is stated the Lemma, by that easily would be solved the class of tasks:

Lemma 1. Let's say that $a, b, n \in N$ and $(a, b) = 1$. Let's prove that if $a^n \equiv 0 \pmod{|a-b|}$, or $b^n \equiv 0 \pmod{|a-b|}$, then $|a - b| = 1$.

Let's solve the equations (I – X) in natural *x*, *y* numbers:

I. $\left(\frac{x+y}{2}\right)^z = x^z - y^z$;

II. $(x+y)^z = (2x)^z + y^z$;

III. $(x+y)^z = (3x)^z + y^z$;

IV. $(y-x)^{x+y} = x^y$, $(y > x)$;

V. $(y-x)^{x+y} = y^x$, $(y > x)$;

VI. $(x+y)^{x-y} = x^y$;

VII. $(x+y)^{x-y} = y^x$;

VIII. $(x+y)^y = (x-y)^x$, $(x > y)$;

IX. $(x-y)^{x+y} = x^{x-y}$;

X. $(x+y)^{x-y} = (x-y)^x$, $(y > x)$.

theorem 5. If a, $b \in N$, $(a, b) = 1$, then each divisor of $(a^2 + ab + b^2)$ will be as similar.

It should be mentioned that specialists of theory of numbers have found a way to solve a number of Diophantine equations and found that many of them cannot be solved, but the general method has not yet been found. By Hilbert were stated the Task on the possibility of the existence of a universal method. In 1970, the famous Russian mathematician I. Matyasevich proved that there could be no such general rule.

In the work also us cist in of our interesting study on Fibonacci numbers. Is introduced the concept of Fibonacci number and is found some of their properties.

**Basic part**

Theorem 1. Let's say that $p$ and $q$ are odd prime numbers and $p = 2q + 1$, from the various prime divisors of number $A = \dfrac{p^q - 1}{p - 1}$ only one of them is less then $p$. $A$ has at least two different divisors.

**Proof**

Let's firstly consider the equation

$$A = q^t, \tag{1}$$

where $A = \dfrac{p^q - 1}{p - 1}$, $p = 2q + 1$, $t \in N$ and $p$, $q$ – are odd prime numbers.

Due introduction we will obtain

$$(2q+1)^q - 1 = 2q^{t+1}. \tag{2}$$

Due application of binomial theorem's formula from (2) we will obtain

$$(2q)^q + C_q^1 (2q)^{q-1} + \ldots + C_q^{q-1} \cdot (2q) = 2q^{t+1}. \tag{3}$$

The left side of (3) would be contained on $2q^2$ and is not divided on than $2q^3$, thus $t + 1 \le 2 \Leftrightarrow t \le 1$. I.e. from (3) we will have

$$(2q+1)^q - 1 > (2q)^q > 2q^2 \ge 2q^{t+1}. \tag{4}$$

From (4) follows the presentation of $A$ as (6.1) is impossible. From the system

$$\begin{cases} p - 1 = 2q \\ \left( \dfrac{p^q - 1}{p - 1}, p - 1 \right) = (p-1, q) = (2q, q) = q \\ A = 1 + p + p^2 + \cdots + p^{q-1} = \dfrac{p^q - 1}{p - 1} > p - 1 = 2q \end{cases} \Rightarrow \begin{cases} A \equiv 0 (\operatorname{mod} q), \\ A > 2q. \end{cases} \tag{5}$$

From the (1) and (5) follow that $A$ has at least two different prime divisors. Let's now say that $A = p_1^{\alpha_1}, p_2^{\alpha_2}, \cdots, p_k^{\alpha_k}$, where $\alpha_1, \alpha_2, \ldots, \alpha_k \in N$, $p_1, p_2, \ldots, p_k -$ are different prime numbers and $p_1 < p_2 < \cdots < p_k$.

Let's say that to $p$ belongs $m$ exponent with module $p_i$, where $p_i \neq q$ and $i$ = 1, 2,..., $k$, then

$$\begin{cases} p^q - 1 = (p-1)(1 + p + p^2 + \cdots + p^{q-1}), \\ 1 + p + p^2 + \cdots + p^{q-1} \equiv 0 (\operatorname{mod} p_i) \\ p^m \equiv 1 (\operatorname{mod} p_i), m - \text{is exponent} \\ q - \text{is odd prime number} \end{cases} \Rightarrow \begin{cases} p^q \equiv 1 (\operatorname{mod} p_i), \\ p^m \equiv 1 (\operatorname{mod} p_i), \\ m - \text{is exponent} \\ q - \text{is odd prime number} \end{cases} \Rightarrow$$

$$\Rightarrow \begin{cases} q \equiv 0 (\operatorname{mod} m), \\ q - \text{is odd prime number} \end{cases} \Rightarrow \left[ \begin{array}{l} m = 1, \\ m = q. \end{array} \right. \qquad (6)$$

As $p_i \neq q$, are odd prime numbers and $p - 1 = 2q \not\equiv 0 (\operatorname{mod} p_i)$, thus $m \neq 1$. I.e. $m = q$.

As exponent $q$ of mutually prime $p$ number at module $p_i$ will divide $\varphi(p_i) = p_i - 1$, so $p_i - 1 \equiv 0 (\operatorname{mod} q) \Rightarrow q < \frac{p_i - 1}{2}$, i.e. $q = p_1$.

$$\begin{cases} q < \dfrac{p_i - 1}{2} \\ q = \dfrac{p-1}{2} \end{cases} \Rightarrow \frac{p_i - 1}{2} > \frac{p-1}{2} \Rightarrow p < p_i; \; i = 1, 2, ..., k.$$

I.e. $A$ has different from $q$ all prime divisors greater than $p$. At the same time $q$ is divisor of $A$. Thus finally we have:

$$q = p_1 < p < p_2 < p_3 < \ldots < p_k. \quad \text{Q.E.D.}$$

Theorem 2. Let's say that $p$ are $q$ odd prime numbers and $p < 2q + 1$, then all prime divisors of $A = \frac{p^q - 1}{p - 1}$ number is greater than $p$.

Let's say that $A = \frac{p^q - 1}{p - 1} = p_1^{\alpha_1}, p_2^{\alpha_2}, \cdots, p_k^{\alpha_k}$, where $p_1, p_2, \cdots, p_k -$ are different odd prime numbers, $p_1 < p_2 < \cdots < p_k$ and $\alpha_1, \alpha_2, \ldots, \alpha_k \in N$. Let's also say that $p$ belongs to exponent $m$ with $p_i$ module $i$ = 1, 2,..., $k$. We have

$$\begin{cases} p^q - 1 = (p-1)(1 + p + p^2 + \cdots + p^{q-1}), \\ 1 + p + p^2 + \cdots + p^{q-1} \equiv 0 (\operatorname{mod} p_i) \\ p^m \equiv 1 (\operatorname{mod} p_i), m - \text{is exponent} \\ q - \text{is odd prime number} \end{cases} \Rightarrow \begin{cases} p^q \equiv 1 (\operatorname{mod} p_i), \\ p^m \equiv 1 (\operatorname{mod} p_i), \\ m - \text{is exponent} \\ q - \text{is odd prime number} \end{cases} \Rightarrow$$

$$\Rightarrow \begin{cases} q \equiv 0 (\operatorname{mod} m), \\ q - \text{is odd prime number} \end{cases} \Rightarrow \left[ \begin{array}{l} m = 1 \\ m = q. \end{array} \right.$$

$$\begin{cases} p < 2q + 1, \\ p, q - \text{are odd prime numbers} \end{cases} \Rightarrow \begin{cases} p - 1 < 2q, \\ p, q - \text{are odd prime numbers} \end{cases}$$

$$\Rightarrow \begin{cases} (p-1, q) = 1, \\ q > \dfrac{p-1}{2}. \end{cases}$$

$$\begin{cases} \left( \dfrac{p^q - 1}{p-1}, p-1) \right) = (p-1, q) = 1, \\ \dfrac{p^q - 1}{p-1} \equiv 0 (\operatorname{mod} p_i). \end{cases} \Rightarrow$$

$\Rightarrow p - 1 \not\equiv 0 (\operatorname{mod} p_i) \Rightarrow p \not\equiv 1 (\operatorname{mod} p_i) \Rightarrow m \neq 1$. I. e. $m = q$.

As exponent $q$ of mutually prime $p$ number at module $p_i$ will divide $\varphi(p_i) = p_i - 1$, thus $p_i - 1 \equiv 0 (\operatorname{mod} q) \Rightarrow q < \dfrac{p_i - 1}{2}$.

We have

$$\begin{cases} q > \dfrac{p-1}{2}, \\ q < \dfrac{p_i - 1}{2}, \\ p_1 < p_2 < \cdots < p_k. \end{cases} \Rightarrow \begin{cases} \dfrac{p-1}{2} < \dfrac{p_i - 1}{2}, \\ p_1 < p_2 < \cdots < p_k. \end{cases} \Rightarrow \begin{cases} p < p_i, \\ p_1 < p_2 < \cdots < p_k. \end{cases} \Rightarrow$$

$\Rightarrow p < p_1 < p_2 < \cdots < p_k$. Q.E.D.

Theorem 3. Let's say that $q$ is odd prime number and $p \in N \setminus \{1\}$, $p \in ]1; q] \cup [q+2; 2q]$, then from the various prime divisors of number $A = \dfrac{p^q - 1}{p-1}$, taken separately, only one is more than $p$.

Proof.

Let's say that $A = \dfrac{p^q - 1}{p-1} = p_1^{\alpha_1}, p_2^{\alpha_2}, \cdots, p_k^{\alpha_k}$, where $p_1, p_2, \cdots, p_k -$ are different odd prime divisors. At the same time, $p_1 < p_2 < \cdots < p_k$ and $\alpha_1, \alpha_2, \ldots, \alpha_k \in N$. Let's alps say that $p$ belongs to exponent $m$ with module $p_1$, then

$$\begin{cases} p^q - 1 = (p-1)(1+p+p^2+\cdots+p^{q-1}), \\ \dfrac{p^q-1}{p-1} \equiv 0 (\bmod p_i), \\ p^m \equiv 1 (\bmod p_1), m - \text{is exponent} \\ q - \text{is odd prime number} \end{cases} \Rightarrow \begin{cases} p^q \equiv 1 (\bmod p_i), \\ p^m \equiv 1 (\bmod p_i), \\ m - \text{is exponent,} \\ q - \text{is odd prime number.} \end{cases} \Rightarrow$$

$$\Rightarrow \begin{cases} q \equiv 0 (\bmod m), \\ q - \text{is odd prime number} \end{cases} \Rightarrow \left[ \begin{array}{l} m = 1, \\ m = q. \end{array} \right.$$

Accordingly of task conditions we have:

$$\begin{cases} p \in N, \\ p \neq q+1, p \neq 2q+1 \\ q - \text{is odd prime number.} \end{cases} \Rightarrow \begin{cases} p \in N, \\ p-1 \neq q, p-1 \neq 2q, \\ q - \text{is odd prime number.} \end{cases} \Rightarrow$$

$$\Rightarrow \begin{cases} (p-1, q) = 1, \\ q - \text{is odd prime number.} \end{cases}$$

$$\begin{cases} \left( \dfrac{p^q-1}{p-1}, p-1 \right) = (p-1, q) = 1, \\ \dfrac{p^q-1}{p-1} \equiv 0 (\bmod p_1). \\ p_1 - \text{is odd prime number.} \end{cases} \Rightarrow p-1 \not\equiv 0 (\bmod p_1) \Rightarrow p \not\equiv 1 (\bmod p_1) \Rightarrow m \neq 1.$$

I.e. $m = q$.

As the exponent $q$ of mutually prime number $p$ at module $p_1$ will be divide $\varphi(p_1) = p_1 - 1$ thus

$$\begin{cases} p_1 - 1 \equiv 0 (\bmod q), \\ p, q - \text{is odd prime number} \end{cases} \Rightarrow q < \frac{p_1 - 1}{2}.$$

I.e. we have:

$$\begin{cases} q > \dfrac{p-1}{2}, \\ q < \dfrac{p_1 - 1}{2}, \\ p_1 < p_2 < \cdots < p_k. \end{cases} \Rightarrow \begin{cases} \dfrac{p-1}{2} < \dfrac{p_1-1}{2}, \\ p_1 < p_2 < \cdots < p_k. \end{cases} \Rightarrow \begin{cases} p < p_i, \\ p_1 < p_2 < \cdots < p_k. \end{cases} \Rightarrow$$

$\Rightarrow p < p_1 < p_2 < \cdots < p_k$. Q.E.D.

Theorem 4. Let's say that $q$ is odd prime number and $p \in \{q + 1; 2q + 1\}$, then from the various prime divisors of number $A = \frac{p^q - 1}{p - 1}$, taken separately, only one is less than $p$. $A$ has at least two different prime divisors.

Proof.

Let's firstly consider the equation

$$A = q^t, \tag{7}$$

where $A = \frac{p^q - 1}{p - 1}$, $p = kq + 1$, $k \in \{1; 2\}$, $t \in N$ and $q$- is odd prime number.

By introduction we will obtain

$$(kq+1)^q - 1 = kq^{t+1},\ k \in \{1; 2\}. \tag{8}$$

Due the formula of binomial theorem from (8) we will obtain

$$(kq)^q + C_q^1 (kq)^{q-1} + \ldots + C_q^{q-1}(kq) = kq^{t+1}. \tag{9}$$

The left side of (9) is divided on $kq^2$ and is not divided on $kq^3$, so $t + 1 \le 2 \Rightarrow t \le 1 \Rightarrow t = 1$. I.e. we will have

$$(kq+1)^q - 1 > (kq)^2 \ge kq^{t+1}. \tag{10}$$

From (10) follows that representation of $A$ as (7) is impossible.

From the system

$$\begin{cases} p - 1 = kq, k \in \{1;2\}, \\ \left( \dfrac{p^q - 1}{p - 1}, p - 1 \right) = (p - 1, q) = (kq, q) = q, \Rightarrow \begin{cases} A \equiv 0 (\bmod q), \\ A > kq. \end{cases} \\ A = \dfrac{p^q - 1}{p - 1} > p - 1 = kq. \end{cases} \tag{11}$$

From the (11) follows that $A$ has at least two different prime divisors. Let's say that $A = p_1^{\alpha_1}, p_2^{\alpha_2}, \cdots, p_k^{\alpha_k}$, where $p_1, p_2, \cdots, p_k -$ are different prime numbers, $p_1 < p_2 < \cdots < p_k$, $\alpha_1, \alpha_2, \ldots, \alpha_k \in N$. Let's also say that $p$ belongs to exponent $m$ with module $p_i$. There $p_i \ne q$, $i = 1, 2, \ldots, k$ then we will have

$$\begin{cases} p^q - 1 = (p-1)(1 + p + p^2 + \cdots + p^{q-1}), \\ 1 + p + p^2 + \ldots + p^{q-1} \equiv 0 (\bmod p_i), \\ p^m \equiv 1 (\bmod p_1), m - \text{is exponent}, \\ q - \text{is odd prime number}, q \ne p_i. \end{cases} \Rightarrow \begin{cases} p^q \equiv 1 (\bmod p_i), \\ p^m \equiv 1 (\bmod p_i), \\ m - \text{is exponent}, \\ q - \text{is odd prime number}. \end{cases} \Rightarrow$$

$$\Rightarrow\begin{cases} q \equiv 0 (\bmod m), \\ q-\text{is odd prime number.}\end{cases} \Rightarrow \left[\begin{array}{l} m=1, \\ m=q. \end{array}\right.$$

$$\begin{cases} \dfrac{p^q-1}{p-1} \equiv 0(\bmod p_i), \\ \left(\dfrac{p^q-1}{p-1}, p-1)\right) = q, \Rightarrow \quad p-1 \not\equiv 0(\bmod p_i) \Rightarrow m \neq 1. \text{ I.e. } m=q. \\ p_i \neq q. \end{cases}$$

As at module $p_i$ the exponent of mutually prime number *p* will be divide $\varphi(p_2) = p_2 - 1$, thus

$$\begin{cases} p_i - 1 \equiv 0(\bmod q), \\ p_i, q - \text{are odd prime numbers.} \end{cases} \Rightarrow q < \frac{p_i - 1}{2}. \text{ I. e. } i \neq 1 \text{ and } q = p_1.$$

$$\begin{cases} q < \dfrac{p_i - 1}{2}, \\ q = \dfrac{p-1}{k}, k \in \{1;2\}, \end{cases} \Rightarrow \frac{p-1}{2} < \frac{p_i - 1}{2} \Rightarrow p < p_i, \ i = 2, 3, ..., k.$$

Finally we have: $q = p_1 < p < p_2 < p_3 < \cdots < p_k$ . Q.E.D.

Lemma 1. Let's say that *a*, *b*, $n \in N$ and (*a*, *b*) = 1. Let's prove that if $a^n \equiv 0 \ (\bmod |a-b|)$, or $b^n \equiv 0 \ (\bmod |a-b|)$, then $|a-b| = 1$.

Proof**.**

$$\begin{cases} a, b \in N, a \neq b, (a,b) = 1 \\ |a-b| \, a^4 \text{ an } |a-b| \, b^4 \end{cases} \Rightarrow \begin{cases} (|a-b|, a^n) = (|a-b|, b^4) = 1 \\ a, b \in N, a \neq b, (a,b) = 1 \end{cases} \Rightarrow |a-b| = 1.$$

Task 1. Let's solve the equation $2^x = \dfrac{y^z - 1}{y - 1}$ in natural numbers. In addition *y* will be the prime number.

Solution:

Let's consider the equation

$$2^x = \frac{y^z - 1}{y - 1}, \tag{12}$$

where *y* is prime number and $x, \ z \in N$.

If *z* = 1, then from (12) $2^x = 1 \Leftrightarrow x = 0 \notin N$ ;

If *z* > 1 and *y* = 2, then the left side of (6.12) is even, and the right side is odd number. I.e. *y* is odd prime number;

If $z > 1$ is odd number and $y$ is odd prime number, then the right side of (12) is odd number, and the keft side – is even that is impossible;

If $z$ is composite number that has the different from 1 odd divisor, i.e. let's say that $z = a \cdot b$, where $a, b \in N$, $b > 1$ and $b$- is odd, then

$$2^x = \frac{(y^a)^b - 1}{y - 1} = \frac{y^a - 1}{y - 1} \cdot (1 + y^a + (y^a)^2 + \ldots + (y^a)^{b-1}). \tag{13}$$

In the brackets of right side of (13) we have odd quantity of odd numbers sum. I.e. this multiplier is odd and in the left side we have $2^x$. This is impossible. Thus (12) has not in this case the solution in natural numbers;

If $z = 2^n$, where $n \in N \setminus \{1\}$, then (12) will be as

$$2^x = \frac{y^{2^n} - 1}{y - 1} = (y+1)(y^2+1)(y^4+1) \cdot \cdots \cdot (y^{2^{n-1}} + 1). \tag{14}$$

In this case (14) has at least two different multipliers ($y$ + 1) and $(y^2 + 1)$, for (12) having the solution in natural numbers, for this ($y$ + 1) and $(y^2 + 1)$ also simultaneously will represent the natural powers of 2.

I. e.

$$\begin{cases} y + 1 = 2^k, \\ y^2 + 1 = 2^l, \end{cases} \tag{15}$$

where $k, l \in N$, $k > 1$, $l > 1$.

From the first equation of (15) $y$ – is odd prime number, thus by dividing $(y^2 + 1)$ on 4 in remainder gives 2, and the second equation in right side of (15) $2^l$ - is dividing without remainder on 4. I.e. in this case (12) has not the solution in natural numbers. It remains to consider the case when $z = 2$.

In this case

$$2^x = \frac{y^2 - 1}{y - 1} = y + 1 \Rightarrow y = 2^x - 1. \tag{16}$$

It would be easily shown that for a prime number $y$ the $x$ must necessarily be prime. I.e. (12) has a solution according to the assignment if and only if $x = P$ – is prime number and $y = P$ – is a Mersenne prime, or $P = 2^p - 1$ – is a prime number, and $z$ is necessarily equal to 2.

From this it would be easily shown that if

$$2^n = \left(\frac{P_1^{\alpha_1} - 1}{P_1 - 1}\right) \cdot \left(\frac{P_2^{\alpha_2} - 1}{P_2 - 1}\right) \cdot \cdots \cdot \left(\frac{P_k^{\alpha_k} - 1}{P_k - 1}\right),$$

where $P_1, P_2, \ldots, P_k$ – are prime numbers and $n$, $\alpha_1$, $\alpha_2, \ldots, \alpha_k \in N$.

Then obligatory:

$$\alpha_1 = \alpha_2 = \ldots = \alpha_k = 2,\ P_1 = 2^{p_1} - 1,\ P_2 = 2^{p_2} - 1, \cdots, P_k = 2^{p_k} - 1,$$

where $p_1, p_2, \ldots, p_k$ – are prime numbers, $p_1 + p_2 + \ldots + p_k = n$, $P_1, P_2, \ldots, P_k$ – are the Mersenne primes,

Task 2. Let's solve the equation $3^x = \dfrac{y^z - 1}{y - 1}$ in natural $x, y, z$ numbers. In addition $y$ – must be prime number.

Solution:

Let's consider the equation

$$3^x = \frac{y^z - 1}{y - 1}, \tag{17}$$

where $y$ – is a prime number and $x, z \in N$.

If $z = 1$, then from (17) $3^x = 1 \Leftrightarrow x = 0 \notin N$.

It is obvious that $y \neq 3$. Let's say that $y$ is odd prime number. Then it is obvious that it is represented as $6a - 1$ or $6a + 1$, where $a \in N$.

If $y = 6a - 1$ and $z$ – is odd, then then at dividing of right side of (17) as remainder will be 1, and the left side – is multiple to 3 that is impossible.

If $y = 6a - 1$ and $z = 2n$, $n \in N$, then from (17) we will have:

$$3^x = \frac{(6a-1)^{2n} - 1}{(6a-1) - 1} = \frac{(6a-1)^2 - 1}{(6a-1) - 1} \cdot ((6a-1)^{2(n-1)} + (6a-1)^{2(n-2)} + \ldots + 1) =$$

$$= 6a((6a-1)^{2(n-1)} + (6a-1)^{2(n-2)} + \ldots + 1),$$

that is impossible, as the left side of this equation is not divided on 2.

If $y = 6a + 1$ and $z$ – are odd, then from (17) we will obtain that $z$ is dividing on 3. I.e. $z=3q$, $q \in N$. By introducing in (17) we will have

$$3^x = \frac{(y^3)^1 - 1}{y - 1} \cdot ((y^3)^{q-1} + (y^3)^{q-2} + \ldots + 1) = 3(3(4a^2 + 2a) + 1)((y^3)^{q-1} + \ldots + 1) \Rightarrow$$

$$\Rightarrow 3^x \equiv 0(\mathrm{mod}(3(4a^2 + 2a) + 1),$$

that is impossible.

If $y = 6a + 1$ and $z$ – are even, let's say that $z = 2n$, $n \in N$, then from (17) we have

$$3^x = \frac{(6a+1)^2 - 1}{(6a+1) - 1} \cdot ((6a+1)^{2(n-1)} + (6a+1)^{2(n-2)} + \ldots + 1) =$$

$$= 2(3a+1)((6a+1)^{2(n-1)}+\ldots+1) \Rightarrow 3^x \equiv 0(\mathrm{mod}\,2),$$

that is impossible;

It remains to consider the case $y = 2$. In this case, from (17) we will obtain

$$3^x = 2^z - 1. \tag{18}$$

Let's say that $x > 1$.

If $x$ –odd, then from (18) we will have:

$$2^z = 3^x + 1 = (3+1)(3^{x-1} - 3^{x-2} + \ldots + 1). \tag{19}$$

In the second bracket of right side of (19) we have odd quantity of odd numbers sum, i.e. odd number, thus (19) has not the sp; in natural numbers;

If $x = 2n$, $n \in N$, then from (18) we will obtain:

$$2^z = 3^{2n} + 1 = (4-1)^{2n} + 1 = 4k + 2, \tag{20}$$

where $k \in N$.

At dividing of right side of (20) on 4 the remainder gives 2, and the left side is dividing on 4, thus also in this case will not the solution of (18):

By us was considered the case for $x = 1$. In this case we will have

$$3^1 = \frac{y^z - 1}{y - 1}. \tag{21}$$

It is obvious that (21) has unique solution $y = z = 2$. Finally we will obtain that the solution of (17), by stated condition is only: $x = 1$; $y = z = 2$.

Task 3. Let's solve the equation

$$P^x = \frac{y^z - 1}{y - 1},$$

where $P \in \{5; 7; 11; ...\}$ – is a prime number, $x, z \in N$ and $y$ – are prime numbers.

Solution:

Let's consider the equation

$$P^x = \frac{y^z - 1}{y - 1}, \tag{22}$$

where $P \in \{5; 7; 11; ...\}$ –is prime number, $y$ – is prime and $x, z \in N$.

Firstly let's show that if $(y^n - 1)$ is dividing on $(y^m - 1)$, where $y \neq 1$ and $m, n \in N$, then obligatory $n$ is dividing on $m$.

Proof.

Let's say that $n = am + r$, where $m, n, a, r \in N$ and $0 \le r < m$. Then we will have

$$y^n - 1 = y^{am} \cdot y^r - 1 = y^r(y^{am} - 1) + (y^r - 1) = y^r(y^m - 1)(y^{m(a-1)} + \cdots + 1) + (y^r - 1). \quad (23)$$

For fulfilment of condition of task, is necessary to divide $(y^r - 1)$ on $(y^m - 1)$, but it is impossible, when $r = 0$ and then we will have $n = am$, where $a, m \in N$. I.e. finally we will obtain provable or $n \equiv 0 (\mathrm{mod}\ n)$. Q.E.D.

From the proved and taking into account that $P$- is a prime number, it is easy to show that (22) has the solution by stated condition, if $z$- is a prime number, or $z = 2^n$, where $n \in N$.

If $z = 2^n$, $n \in N$, then (22) will be as

$$P^x = \frac{y^{2^n} - 1}{y - 1}. \quad (24)$$

For $n = 1$ $P^x = \frac{y^2 - 1}{y - 1} = y + 1 \Rightarrow y = P^x - 1$ if we take into account that $y$ and $P$- are prime numbers, then $P = 2$ and $x = P_1$, where $P_1 -$ is a prime number and $y = 2^{P_1} - 1$ – is Messene prime;

If $n > 1$, then (22) will be as

$$P^x = \frac{y^{2^n} - 1}{y - 1} = (1 + y)(1 + y^2) \cdot \cdots \cdot (1 + y^{2^{n-1}}). \quad (25)$$

In this case in the right side (25) has at least two multipliers $(y + 1)$ and $(y^2 + 1)$. For having to (25) the solution in natural numbers, for this either $(y + 1)$ and $(y^2 + 1)$ simultaneously will be presented by any natural power of $P$. I.e.

$$\begin{cases} y + 1 = P^a, \\ y^2 + 1 = P^b \end{cases} \Rightarrow P^b = (y + 1)^2 - 2y = P^{2a} - 2y. \quad (26)$$

As $(P, y) = 1$, thus from (26) follows that $P \equiv 0 (\mathrm{mod}\ 2)$. But $P$- is prime number, so accordingly of task $P = 2$ and (26) would be written down as

$$2^b = 2^{2a} - 2y \Rightarrow y = 2^{b-1}(2^{2a-b} - 1). \quad (27)$$

$y$- is a prime number, thus from the (27) we will have: $b = 1$, $2a - 1 = P_1$ – are prime numbers and $y = 2^{P_1} - 1$ – represents a Messene prime.

It would be mentioned that

$$y < P^{\frac{x}{z-1}} < y + 1, \quad (28)$$

Therefore if $\frac{x}{z-1} \in N$, then (22) has not the solution!

Let's now the case, when $z$ – is a prime number.

I.e. we will have

$$P^x = \frac{y^z - 1}{y-1},$$

where $x \in N$ and $P$, $z$, $y$ – are prime numbers.

We have

$$\begin{cases} \left(\dfrac{y^z-1}{y-1}, y-1\right) = (y-1, z). \\ y, z, P - \text{prime numbers}, \quad \Rightarrow P-1 \equiv 0 (\operatorname{mod} z) \Rightarrow P \equiv 1 (\operatorname{mod} z) \Rightarrow P = 2kz+1, \\ \dfrac{y^z-1}{y-1} \equiv 0 (\operatorname{mod} P), \end{cases}$$

where $P$, $z$ – are prime numbers and $k \in N$.

Let's consider the several cases

**1** If $P = 5 = 2 \cdot 2 + 1$, i.e. $z = 2$, then from (22) we will have

$$5^x = \frac{y^2-1}{y-1} = y+1 \Rightarrow y = 5^x - 1 \equiv 0 (\operatorname{mod} 4).$$

This is impossible, as $y$ – is a prime number.

**2** If $P = 7 = 2 \cdot 3 + 1$, i.e. $z = 3$. Then from (22) we have

$$7^x = \frac{y^3-1}{y-1} = y^2 + y + 1. \tag{29}$$

Let's consider the identities:

$$(a^2+ab+b^2)(c^2+cd+d^2) = (ac-bd)^2 + (ac-bd)(ad+bc+bd) +$$
$$+(ad+bc+bd)^2 = (ad-bc)^2 + (ad-bc)(ac+bd+bc) + (ac+bd+bc)^2; \tag{30}$$

$$(a^2+ab+b^2)^2 = (a^2-b^2)^1 + (b^2+2ab)(a^2-b^2) + (b^2+2ab)^2 =$$
$$= (b^2-a^2)^2 + (b^2-a^2)(a^2+2ab) + (a^2+2ab)^2; \tag{31}$$

$$(a^2+ab+b^2)^3 = (a^3-3b^2a-b^3)^2 + (a^3-3b^2a-b^3)(3ab(a+b)) + (3ab(a+b))^2. \tag{32}$$

$$7^2 = 1^2 + 1 \cdot 2 + 2^2.$$

By this formulae we will have:

$$(1^2+1\cdot 2+2^2)^2 = 3^2+3\cdot 5+5^2, \quad (1^2+1\cdot 2+2^2)^3 = 1^2+1\cdot 18+18^2.$$

$$(1^2+1\cdot 2+2^2)^5 = (1^2+1\cdot 18+18^2)(3^2+3\cdot 5+5^2) = 7^2+7\cdot 126+126^2.$$

Let's say that we have:

$$\begin{cases}1+p+p^2=(1+2+2^2)(a^2+ab+b^2)\\ a^2+ab+b^2=(1+2+2^2)(c^2+cd+d^2)\\ a,b,c,d\in N,\ b>a,\ c>d\end{cases}\Rightarrow\begin{cases}|1\cdot b-2\cdot a|=1\\ a\cdot 1+2\cdot b+2\cdot a=p\\ |1\cdot c-2\cdot d|=a\\ 1\cdot d+2\cdot c+2\cdot d=b\\ a,b,c,d\in N,\ b>a,\ c>d\end{cases}$$

From this we obtain: $\begin{cases}b=2a\pm 1, p=7a\pm 2\\ c=2d+a\\ 7d=\pm 1\\ a,b,c,d\in N,\ b>a,\ c>d\end{cases}$ or $\begin{cases}b=2a\pm 1, p=7a\pm 2\\ c=2d-a>d\\ 7d=4a\pm 1<4d+1\\ a,b,c,d\in N,\ b>a,\end{cases}$, from that no one has the solution in natural numbers. I.e. $\dfrac{1+p+p^2}{(1^2+1\cdot 2+2^2)^3}\notin N$. Thus last presentation of $(1+p+p^2)$ as $(1^2+1\cdot2+2^2)^m$ for $m=3$ $(1^2+1\cdot2+2^2)^3 =\ 1^2 + 1\cdot18 + 18^2$.

I.e. the solutions of (29) are only: $x = 1$; $y = 2$; $z = 3$.

Task 4. Let’s solve the equation

$$\left(\frac{x+y}{2}\right)^z=x^z-y^z \tag{33}$$

in natural numbers.

Solution:

Let’s rewrite (33) as

$$(x+y)^z=2^z(x^z-y^z)\,. \tag{34}$$

For $z > 1$, from (6.34) we will obtain $x + y = 2x - 2y \Leftrightarrow \mathbf{x} = 3y$. I.e. the solutions of (33) are: $x = 3t$, $y = t$, where $t \in N$;

for $z>1$ let’s say that sols of (34) are: $x_0$, $y_0$, $z_1$ and $x_0=x_1d$, $y_0=y_1d$, where $d=(x_0,y_0)$.

By introduction in (34) and dividing of two sides on $d^{z_1}$ we will obtain

$$(x_1+y_1)^{z_1}=2^{z_1}(x_1^{z_1}-y_1^{z_1})\,. \tag{35}$$

The right side of (25) is a multiple of 2. Therefore $x_1$ and $y_1$ numbers are odd.

I.e.

$$x_1=2a-1 \text{ end } y_1=2b-1, \text{ where } a, b \in N. \tag{36}$$

By introduction of (36) in (35) and simplification we will obtain

$$(a+b-1)^{z_1}=(2a-1)^{z_1}-(2b-1)^{z_1}\,. \tag{37}$$

$$(2a-1)^{z_1}-(2b-1)^{z_1}=((2a-1)-(2b-1))(2a-1)^{z_1-1}+\ldots+(2b-1)^{z_1-1})=$$

$$= 2(a-b)((2a-1)^{z_1-1} + \ldots + (2b-1)^{z_1-1}), \qquad (38)$$

$$(2b-1)^{z_1} = (2a-1)^{z_1} - (a+b-1)^{z_1} = ((2a-1)-(a+b-1))((2a-1)^{z_1-1} + \ldots$$

$$\ldots + (a+b-1)^{z_1-1}) = (a-b)\left((2a-1)^{z_1-1} + \ldots + (a+b-1)^{z_1-1}\right) \qquad (39)$$

From (37), (38) and (39) we have that $(a-b)\setminus(2a-1)^{z_1}$ and $(a-b)\setminus(2b-1)^{z_1}$, but ((2*a* – 1), (2*b* – 1)) = 1, therefore accordingly of Lemma 1 by introduction in *a* – *b* = 1 or $a = b+1$ (37).

$$(2b)^{z_1} = (2b+1)^{z_1} - (2b-1)^{z_1} \Leftrightarrow (2b-1)^{z_1} + (2b)^{z_1} = (2b+1)^{z_1}. \qquad (40)$$

If $z_1 = 1$, then from (40)

$$(2b)^{z_1} = \left((2b+1)-(2b-1)\right)\left((2b+1)^{z_1-1} + \ldots + (2b-1)^{z_1-1}\right) =$$

$$= 2\left((2b+1)^{z_1-1} + \ldots + (2b-1)^{z_1-1}\right) \Rightarrow 2^{z_1-1} \cdot b^{z_1} = (2b+1)^{z_1-1} + \ldots + (2b-1)^{z_1-1}. \qquad (41)$$

If $z_1 > 1$ is an odd natural number, then the right side of (41) is an odd amount of the sum of odd numbers, or an odd number, and the left side is even, so for odd $z_1$ numbers other than 1, this equation has no solution in natural numbers *x*, *y*, *z*.

If $z_1$ is even, i.e. $z_1 = 2t$, then (40) will be as

$$(2b-1)^{2t} + (2b)^{2t} = (2b+1)^{2t}. \qquad (42)$$

If *t* = 1, then we have $(2b-1)^2 + (2b)^2 = (2b+1)^2$, from that $z_1 = b = 2$.

I.e. the solution of (35) are:

$$\begin{cases} x_1 = 3, \\ y_1 = 1, \\ z_1 = 1. \end{cases} \quad \text{or} \quad \begin{cases} x_1 = 5, \\ y_1 = 3, \\ z_1 = 2. \end{cases}$$

Thus the solutions of (33) will be:

$$\begin{cases} x = 3k, \\ y = k, \\ z = 1. \end{cases} \quad \text{or} \quad \begin{cases} x = 5k, \\ y = 3k, \\ z = 2, \end{cases} \qquad (43)$$

where $k \in N$.

If *t* > 1, then from the binomial theorem formula from (42) is easy to see that $b \setminus t$.

By the dividing of both sides of (42) on $(2b)^{2t}$ we will obtain:

$$2 > \left(1 - \frac{1}{2b}\right)^{2t} + 1 = \left(1 + \frac{1}{2b}\right)^{2t} > 1 + C_{2t}^{1}\left(\frac{1}{2b}\right) = 1 + \frac{2t}{2b} = 1 + \frac{t}{b} \Rightarrow 1 > \frac{t}{b} \Rightarrow t < b,$$

but $b \setminus t$, that is impossible.

Finally we have that all sons of (33) gives (43).

Task 5. Let’s solve the equation

$$(x+y)^z=(2x)^z+y^z, \qquad (44)$$

in natural numbers.

Solution:

Firstly let's say that $z > 2$ and solutions of (44) are: $x_0$, $y_0$, $z_1 \in N$, $z_1 > 2$. At the same time $x_0 = x_1 d$, $y_0 = y_1 d$, where $d = (x_0; y_0)$. By dividing of both sides of (44) on $d^{z_1}$ we will obtain:

$$(x_1+y_1)^{z_1}=(2x_1)^{z_1}+y_1^{z_1} \Leftrightarrow y^{z_1}=(x_1+y_1)^{z_1}-(2x_1)^{z_1}=$$

$$=\left((x_1+y_1)-2x_1\right)\left((x_1+y_1)^{z_1-1}+(x_1+y_1)^{z_1-2}(2x_1)+\ldots+(2x_1)^{z_1-1}\right)=$$

$$=(y_1-x_1)\left((x_1+y_1)^{z_1-1}+(x_1+y_1)^{z_1-2}\cdot(2x_1)+\ldots+(2x_1)^{z_1-2}\right)\Rightarrow(y_1-x_1)\setminus y_1^{z_1},$$

In addition, $(x_1, y_1) = 1$. Thus accordingly of Lemma 6.1 $y_1 - x_1 = 1$ or $y_1 = x_1 + 1$.

By introduction in (44) we have

$$(2x_1+1)^{z_1}=(2x_1)^{z_1}+(x_1+1)^{z_1}. \qquad (45)$$

The left side of (45) is odd, there fire $(x_1+1)$ will be odd, i.e. $x_1 = 2k$, $k \in N$. By introduction in (45) we will obtain

$$(4k+1)^{z_1}=(4k)^{z_1}+(2k+1)^{z_1}. \qquad (46)$$

If $z_1$ is odd, then from (46)

$$(4k+1)^{z_1}=(4k+(2k+1))\left((4k)^{z_1-1}-(4k)^{z_1-2}(2k+1)+\ldots+(2k+1)^{z_1-1}\right)\Rightarrow$$

$$\Rightarrow(6k+1)\setminus(4k+1)^{z_1},$$

but $(6k + 1; 4k + 1) = 1$ (indeed $6k + 1 = (4k + 1) + 2k$)). I.e. $z_1$ is even or $z_1 = 2n$, where $n \in N$, $n \neq 1$. Due taking it into account (46) will be as:

$$\left((4k+1)^n\right)^2=\left((4k)^n\right)^2+\left((2k+1)^n\right)^2. \qquad (47)$$

$(4k + 1)$, $4k$ and $(2k + 1)$ pairwise are coprime numbers. Thus would be found such $p, q \in N$, $p > q$ and $(p, q) = 1$, that:

$$(4k+1)^n=p^2+q^2,\ (4k)^n=2pq,\ (2k+1)^n=p^2-q^2. \qquad (48)$$

From (48) we have

$$2q^2=(p^2+q^2)-(p^2-q^2)=(4k+1)^n-(2k+1)^n=$$

$$=\left((4k+1)-(2k+1)\right)\left((4k+1)^{n-1}+\ldots+(2k+1)^{n-1}\right)=$$

$$=2k\left((4k+1)^{n-1}+\ldots+(2k+1)^{n-1}\right)\Rightarrow k\setminus q^2. \qquad (49)$$

As $(p, q) = 1$ and $k \setminus q^2$, therefore $(k, p) = 1$.

$$(4k)^n = 2pq \Rightarrow 2^{2n-1}k^n = pq. \quad (60)$$

If $k = 1$, then (46) will be as:

$$5^{z_1} = 4^{z_1} + 3^{z_1} \quad \text{or}$$

$$\begin{cases} \left(\frac{4}{5}\right)^{z_1} + \left(\frac{3}{5}\right)^{z_1} = 1, \\ z_1 < 2 \Rightarrow \left(\frac{4}{5}\right)^{z_1} + \left(\frac{3}{5}\right)^{z_1} > \left(\frac{4}{5}\right)^2 + \left(\frac{3}{5}\right)^2 = 1, \\ z_1 = 2 \Rightarrow \left(\frac{4}{5}\right)^{z_1} + \left(\frac{3}{5}\right)^{z_1} = \left(\frac{4}{5}\right)^2 + \left(\frac{3}{5}\right)^2 = 1 \\ z_1 > 2 \Rightarrow \left(\frac{4}{5}\right)^{z_1} + \left(\frac{3}{5}\right)^{z_1} < \left(\frac{4}{5}\right)^2 + \left(\frac{3}{5}\right)^2 = 1. \end{cases} \Rightarrow z_1 = 2 \quad (51)$$

If $k > 1$, then from (49) and (50) $k$ is odd and $p = 2^{2n-1}$, $q = k^n$. From (48) we will obtain

$$(4k+1)^n = (2^{2n-1})^2 + (k^n)^2 .$$

From that

$$2^{4n-2} = (4k+1)^n - (k^2)^n = (4k+1-k^2)\left((4k+1)^{n-1} + (4k+1)^{n-2}k^2 + \ldots + (k^2)^{n-1}\right). \quad (52)$$

From (52) follow that $4k+1-k^2 \geq 1 \Leftrightarrow k \in ]0;4[$. At the same time $k$ is other than 1 odd number. I.e. $k = 3$. By introduction in (47) we will obtain

$$13^{2n} = 12^{2n} + 7^{2n}. \quad (53)$$

If $n = 1$, then $13^{2n} = 13^2 \neq 12^2 + 7^2 = 12^{2n} + 7^{2n}$.

If $n > 1$, then from (53) we have

$$12^{2n} = (13^2)^n - (7^2)^n = 169^n - 49^n = (169-49)(169^{n-1} + 160^{n-2}49 + \ldots + 49^{n-1}) =$$

$$= 120(169^{n-1} + \ldots + 49^{n-1}). \quad (54)$$

The right side of (54) is multiple to 5, and the left side not, thus (54) has no solution in natural numbers.

Finally we have that if $z > 2$, then (44) has no solution in natural numbers.

If $z = 1$, then from (44) we will obtain $x + y = 2x + y \Leftrightarrow x = 0$, that is impossible, as $x \in N$.

If $z = 2$, then

$$(x+y)^2 = (2x)^2 + y^2 \Leftrightarrow 3x = 2y, \quad (55)$$

From that we will obtain that all natural solutions of (44) are:

$$x = 2t, \quad y = 3t, \quad z = 2,$$

where $t \in N$.

Task 6. Let's solve the equation

$$(x+y)^z = (3x)^z + y^z\,, \tag{56}$$

in natural numbers.

Solution:

Let's say that the solutions of (56) are

$$x_1, y_1, z_0 \in N \text{ and } (x_1, y_1) = d\,,\; x_1 = x_0 d\,,\; y_1 = y_0 d\,. \tag{57}$$

By introduction of (57) in (56) and simplification we will obtain:

$$(x_0 + y_0)^{z_0} = (3x_0)^{z_0} + y_0^{z_0}\,. \tag{58}$$

From (58) if $z_0 > 1$ we will have:

$$y_0^{z_0} = (y_0 - 2x_0)\left((x_0 + y_0)^{z_0 - 1} + \ldots + (3x_0)^{z_0 - 1}\right). \tag{59}$$

From that

$$y_0^{z_0} \equiv 0(\operatorname{mod}(y_0 - 2x_0))\,. \tag{60}$$

If $y_0$ - is odd, then accordingly of Lemma 1 from (60) we have:

$$y_0 - 2x_0 = 1 \text{ or } y_0 = 2x_0 + 1\,. \tag{61}$$

By introduction of (61) in (58) we will obtain:

$$(3x_0 + 1)^{z_0} = (3x_0)^{z_0} + (2x_0 + 1)^{z_0}\,. \tag{62}$$

If $z_0$ - is odd, then from (62) we will have:

$$(3x_0 + 1)^{z_0} = (5x_0 + 1)\left((3x_0)^{z_0 - 1} - \ldots + (2x_0 + 1)^{z_0 - 1}\right), \tag{63}$$

$$5x_0 + 1 = \left((3x_0 + 1)^{z_0}, 5x_0 + 1\right) = \left((5x_0 + 1 - 2x_0)^{z_0}, 5x_0 + 1\right) =$$

$$= \left((2x_0)^{z_0}, 5x_0 + 1\right) = (2^{z_0}, 5x_0 + 1)\,. \tag{64}$$

I.e. if $z_0$ - is odd $(z_0 \neq 1)$, then for (62) to have a solution, then $x_0$ also will be odd. From (63) and (64) follows that $5x_0 + 1 = 2^{z_0}$, as in the second bracket on the right side of (63) would have an odd amount of sums of odd numbers, or odd number. At the same time $(3x_0 + 1)^{z_0} = (3(2k+1)+1)^{z_0} = 2^{z_0}(3k+2)^{z_0}$ and in the right side even is only $(5x_0 + 1)$ multiplier.

I.e.

$$5x_0 = 2^{z_0} - 1\,. \tag{65}$$

For odd $z_0$ the right side of (65) is ended by 1 or 7, thus (65) will not be executed for odd $x_0$ and $z_0$. I.e. $z_0$ - is an even.

As the equation $(a^n)^4+(b^n)^4=(c^n)^4$ has no solution in $a, b, c, n\in N$ (this know fact is proven due infinite assumptions method), thus even $z_0$ would be as

$$z_0 = 4(n-1)+2,\ n \in N. \tag{66}$$

From (62) and (66) we will have

$$(9x_0^2+6x_0+1)^{2n-1}=(9x_0^2)^{2n-1}+(4x_0^2+4x_0+1)^{2n-1}. \tag{67}$$

If $n > 1$, $n \in N$, then from (67) we will have:

$$(6x_0+1)(9x_0^2+6x_0+1)^{2n-2}+(9x_0^2+6x_0+1)^{2n-3}\cdot(9x_0^2)+\cdots+(9x_0^2)^{2n-2})=$$

$$=(2x_0+1)^{4n-2}\ . \tag{68}$$

As $(6x_0+1, 2x_0+1)=(6x_0+1, 6x_0+2)=(2,6x_0+3)=1$, thus (68) has no solution accordingly of above mentioned condition.

If $z_0=1$, then from (68)

$$(x_0+y_0)=(3x_0)+y_0 \Rightarrow x_0=0\notin N.$$

If $z_0=2$, then from (62)

$$(3x_0+1)^2=(3x_0)^2+(2x_0+1)^2 \Leftrightarrow 9x_0^2+6x_0+1=$$

$$=9x_0^2+4x_0^2+x_0+1 \Leftrightarrow x_0=0\notin N\ .$$

Finally we have that (56) accordingly of this condition has no solution.

Task 7. Let’s solve the equation

$$(y-x)^{x+y}=x^y,\ (y>x) \tag{69}$$

in natural numbers.

Solution:

Let’s say that solutions of stated equation are $x_0,\ y_0$ and at the same time $x_0=x_1d$, $y_0=y_1d$, where $d=(x_0,y_0)$. Due introduction in the stated equation, we will obtain:

$$d^{x_1d}(y_1-x_1)^{x_1d}=\left(\frac{x_1}{y_1-x_1}\right)^{y_1d}. \tag{70}$$

From this

$$d^{x_1}(y_1-x_1)^{x_1}=\left(\frac{x_1}{y_1-x_1}\right)^{y_1}. \tag{71}$$

From (71) $(y_1-x_1)\backslash x_1$. Accordingly of Lemma 1

$$y_1-x_1=1\Rightarrow y_1=x_1+1. \tag{72}$$

By introduction of (72) in (71) and simplification we will obtain $d^{x_1}=x_1^{x_1+1}$ or

$$d = x_1 \cdot x_1^{\frac{1}{x_1}}. \tag{73}$$

It is easy to show that if $n^{\frac{1}{n}}$ – is a rational number, then it obligatory will be natural.

Let's now say that.

$$k,\ n\ \in N,\ n \neq 1 \text{ and } n^{\frac{1}{n}} = k. \tag{74}$$

From (74) $n = k^n \geq 2^n$. But the last is not valid. I.e. $x_1^{\frac{1}{x_1}}$ – is an irrational if $x_1 \in N$, $x_1 \neq 1$.

For $x_1 = 1$ we will assume that (71) has unique pair of solutions:

$$\begin{cases} x_1 = 1 \\ y_1 = 2. \end{cases}$$

290 $\Rightarrow$ $d$ =1 $\Rightarrow$therefore the solutions of (69) are (1, 2).

Task 8. Let's solve the equation

$$(y - x)^{x+y} = y^x \quad (y > x) \tag{75}$$

in natural numbers.

Solution:

Let's say that solutions of (75) are $x_0, y_0$ and at the same time

$$x_0 = x_1 d\,,\ \ y_0 = y_1 d\,, \tag{76}$$

where $d_0 = (x_0, y_0)$.

By introduction of (76) in (75) we will obtain:

$$d^{x_1 d + y_1 d} (y_1 - x_1)^{x_1 d + y_1 d} = d^{x_1 d} \cdot y_1^{x_1 d}. \tag{77}$$

From (77) we have:

$$d^{y_1} (y_1 - x_1)^{x_1} = \left( \frac{y_1}{y_1 - x_1} \right)^{x_1}. \tag{78}$$

The left side of (78) is natural number, thus accordingly coif Lemma 1

$$y_1 = x_1 + 1, \text{ from that } \ y_1 = x_1 + 1. \tag{79}$$

By introduction of (79) in (6.78) we will obtain:

$$d^{x_1 + 1} = (x_1 + 1)^{x_1} \Rightarrow d = \left( \frac{x_1 + 1}{d} \right)^{x_1}. \tag{80}$$

In order for (80) to have a solution, it is necessary to find $k \in N \setminus \{1\}$, for that $x_1 + 1 = kd$. By taking into account (80) we obtain

$$d = k^{kd-1}. \quad (81)$$

But as for each $n \in N$ $n < 2^{2n-1}$, thus (80) has no solution and hence also (75) will have the solution in natural numbers.

Task 9. Let's solve the equation

$$(x+y)^{x-y} = x^y \quad (82)$$

in natural numbers.

solution:

Let's say that solutions of (82) are $x_0, y_0$ and at the same time

$$x_0 = x_1 d, \quad y_0 = y_1 d, \quad (83)$$

where $d = (x_0, y_0)$.

By introduction of (83) in (82) we will obtain:

$$d^{(x_1-y_1)d}(x_1+y_1)^{(x_1-y_1)d} = x_1^{y_1 d} d^{y_1 d}. \quad (84)$$

From (84) we will obtain:

$$\frac{(x_1+y_1)^{x_1-y_1}}{x_1^{y_1}} = d^{2y_1-x_1}. \quad (85)$$

Let's say that $x_1 \neq 1$. Then we will have:

$$(x_1, y_1) = 1 \Rightarrow (x_1+y_1, x_1) = 1 \Rightarrow \left((x_1+y_1)^{x_1-y_1}, x_1^{y_1}\right) = 1. \quad (86)$$

The right side of (85) or its reciprocal number is natural. From (6.86) follows that neither left side of (85) and nor its reciprocal value is possible to be natural number. I.e. for $x_1 \neq 1$ also (85) has no solution in natural numbers.

If $x_1 = 1$, then from (85) easy will be obtained $x_1 = y_1 = d = 1$.

I.e. (82) has the unique solution in natural numbers: $x_0 = y_0 = 1$.

Task 10. Let's solve the equation

$$(x+y)^{x-y} = y^x \quad (87)$$

in natural numbers.

solution:

Let's say that solutions of (87) are $x_0$, $y_0$ and at the same time

$$x_0 = x_1 d, \quad y_0 = y_1 d, \quad (88)$$

where $d = (x_0, y_0)$.

By introduction of (88) in (87) we will obtain:

$$d^{(x_1-y_1)d}(x_1+y_1)^{(x_1-y_1)d} = y_1^{x_1 d}d^{x_1 d}\,, \qquad (89)$$

From that we will obtain:

$$d^{x_1-y_1}(x_1+y_1)^{x_1-y_1} = y_1^{x_1}d^{x_1}\,. \qquad (90)$$

From (90) we will obtain

$$d^{y_1} = \frac{(x_1+y_1)^{x_1-y_1}}{y_1^{x_1}}\,. \qquad (91)$$

$$(x_1,y_1)=1 \Rightarrow (x_1+y_1,y_1)=1 \Rightarrow \left((x_1+y_1)^{x_1-y_1}, y_1^{y_1}\right)=1\,. \qquad (92)$$

The left side of (91) is a natural number, as well as right side, to be natural, from (92) follows that $y_1$ necessary will be equal to 1. I.e. $y_1 = 1$.

By introduction of $d = (x_1+1)^{x_1-1}$ in (91) we finally have that all pairs of solution of (87) will be from

$$\begin{cases} x = t(t+1)^{t-1}, \\ y = (t+1)^{t-1} \\ t \in N \end{cases}$$

Task 11. Let's solve the equation

$$(x+y)^y = (x-y)^x\,,\ (x>y) \qquad (93)$$

in natural numbers.

Solution:

Let's say that solution of (93) is $x_0$ and $y_0$, where $(x_0,y_0)=d$, $x_0 = x_1 d$ and $y_0 = y_1 d$. I.e. $(x_1,y_1)=1$, $x_1,y_1 \in N$.

By introduction in (93) and simplification

$$(x_1+y_1)^{y_1} = (x_1-y_1)^{x_1}d^{x_1-y_1}\,. \qquad (94)$$

I.e.

$$\left(\frac{x_1+y_1}{x_1-y_1}\right) \in N_1 \ \text{ or } \ \left(1+\frac{2y_1}{x_1-y_1}\right) \in N \Rightarrow \frac{2y_1}{x_1-y_1} \in N\ .$$

From the last accordingly of Lemma 1,

$$x_1 - y_1 = 1, \quad (95) \qquad \text{or} \qquad x_1 - y_1 = 2\,. \quad (96)$$

Let's take into account the case (95)

$$x_1 - y_1 = 1 \Rightarrow y_1 = x_1 - 1\ .$$

It is obvious that $x_1 > 1$.

By introduction in (94) and simplification we will obtain:

$$d = (2x_1 - 1)^{x_1 - 1}.$$

I.e. the solutions of (93) are

$$x_0 = x_1(2x_1 - 1)^{x_1 - 1} \text{ and } y_0 = (x_1 - 1)(2x_1 - 1)^{x_1 - 1}, \qquad (97)$$

where $x_1 \in N \setminus \{1\}$.

Let's take into account the case (96)

$$x_1 - y_1 = 2 \Rightarrow y_1 = x_1 - 2.$$

It is obvious that $x_1 \in N \setminus \{1\}$.

By introduction in (94) and simplification, we will obtain

$$(x_1 - 1)^{x_1 - 1} = (2d)^2. \qquad (98)$$

From (98) is clear that $x_1$ is impossible to be even number, thus $x_1 - 1 = (2n)^2$, where $n \in N$, or $x_1 = 4n^2 + 1$.

From (98), $d = \dfrac{(2n)^{4n^2 - 1}}{2}$; i.e. in this case:

$$x_0 = (4n^2 + 1) \cdot \frac{(2n)^{4n^2 - 1}}{2}; \quad y_0 = (4n^2 - 1) \cdot \frac{(2n)^{4n^2 - 1}}{2}. \qquad (99)$$

Finally we have that all solutions of (93) are stated by (97) and (99) formula.

Task 12. Let's solve the equation

$$(x - y)^{x + y} = x^{x - y} \quad (x > y) \qquad (100)$$

in natural numbers.

Solution:

Let's say that solution of equation (100) are $x_0$ and $y_0$, $x_0 = x_1 d$, $y_0 = y_1 d$, where $d = (x_0, y_0)$, $x_1, y_1, d \in N$. By introduction in (100) we will obtain

$$d^{(x_1 + y_1)} \cdot (x_1 - y_1)^{(x_1 + y_1)d} = x_1^{(x_1 - y_1)d} \cdot d^{(x_1 + y_1)d},$$

From that

$$d^{2y_1}(x_1 - y_1)^{2y_1} = \left( \frac{x_1}{x_1 - y_1} \right)^{x_1 - y_1}. \qquad (101)$$

The left side of (101) is a natural number. Thus natural will be also $\left( \dfrac{x_1}{x_1 - y_1} \right)$. Accordingly of this Lemma 1 $x_1 - y_1 = 1 \Leftrightarrow y_1 = x_1 - 1$. By introduction in (101) and simplification we will obtain

$$d^{2(x_1-1)} = x_1 . \tag{102}$$

By application of method of mathematical induction, easy is to prove that if $x_1 \in N \setminus \{1\}$, then $2^{2(x_1-1)} > x_1$, thus from (102) will be obtained that $d$ = 1 and $x_1 = 1$, that is impossible as $x_1, y_1 \in N$ and $x_1 > y_1$. I.e. (100) has no solution in natural numbers.

Task 13. Let's solve the equation

$$(x+y)^{x-y} = (x-y)^x \quad (x > y) \tag{103}$$

in natural numbers

Solution:

Let's say that solutions of (103) are $x_0$ and $y_0$, $x_0 = x_1 d$, $y_0 = y_1 d$, where $x_1, y_1, d \in N$; $(x_0, y_0) = d$. By introduction in (103) and simplification we will obtain:

$$(d(x_1+y_1))^{(x_1-y_1)d} = (d(x_1-y_1))^{x_1 d}, \tag{104}$$

from that we obtain

$$d^{x_1-y_1}(x_1+y_1)^{x_1-y_1} = d^{x_1}(x_1-y_1)^{x_1 d}, \tag{105}$$

Hence we have

$$\left(\frac{x_1+y_1}{x_1-y_1}\right)^{x_1-y_1} = d^{y_1}(x_1-y_1)^{y_1} . \tag{106}$$

I.e. $\frac{x_1+y_1}{x_1-y_1} = 1 + \frac{2y_1}{x_1-y_1} \in N$, or accordingly of Lemma 1 we have

$$\left[\begin{array}{l} x_1 - y_1 = 1, \\ x_1 - y_1 = 2. \end{array}\right. \Leftrightarrow \left[\begin{array}{l} y_1 = x_1 - 1, \quad (107) \\ y_1 = x_1 - 2. \quad (108) \end{array}\right.$$

By introduction of (107) in (108) we will obtain $d^{x_1-1} = 2x_1 - 1$, from that $d = 1$ and $x_1 = 1$, as in contrary

$$d^{x_1-1} \geq 2^{x_1-1} \geq 2x_1 - 1. \tag{109}$$

Due the mathematical induction method would is easy to be proven the inequality

$$2^{x_1-1} > 2x_1 - 1, \text{ where } x_1 \in N \setminus \{1\}.$$

I.e. in this case (103) has no solution because $y_1 = 1 - 1 = 0 \notin N$.

By introduction of (108) in (106) we will obtain

$$\left(\frac{2x_1-2}{2}\right)^2 = d^{x_1-2} \cdot 2^{x_1-2} \Leftrightarrow (x_1-1)^2 = (2d)^{x_1-2} . \tag{110}$$

If $\begin{cases} d > 1, \\ x_1 > 2, \end{cases}$ then by the method of mathematical induction it is easy to prove that $(2d)^{x_1-2} > (x_1-1)^2$.

If $d=1$, then from (110)

$$2^{x_1-2}=(x_1-1)^2 . \qquad (111)$$

From (111) either $(x_1-1)$ and $(x_1-2)$ simultaneously will be even, that is impossible. I.e. (111) has no solution in natural numbers and therefore, (103) also will have the solution in natural numbers.

Theorem 5. If $(a, b) = 1$, then $a^2+ab+b^2$ and each divisor of (112) will be as (112).

Firstly let's prove the Lemma 2 and Lemma 3.

Lemma 2. If $(a, b) = 1$ and $a^2+ab+b^2$ type number is divided on $c^2+cd+d^2=p$ ($a, b, c, d \in N$) prime number, then the divider also would be present as following type of number.

Proof.

Let's consider the identities:

$$\begin{cases}(a^2+ab+b^2)(c^2+cd+d^2)=(ac-bd)^2+(ac-bd)(ad+bc+bd)+(ad+bc+bd)^2; & (113)\\ (a^2+ab+b^2)(c^2+cd+d^2)=(ad-bc)^2+(ad-bc)(ac+bd+bc)+(ac+bd+bc)^2; & (114)\\ (ac-bd)(ac+bd+bc)=c^2(a^2+ab+b^2)-b^2(c^2+cd+d^2), & (115)\\ (ad-bc)(ad+bc+bd)=d^2(a^2+ab+b^2)-b^2(c^2+cd+d^2). & (116)\end{cases}$$

Let's say that $a^2+ab+b^2\equiv 0(\mathrm{mod}(c^2+cd+d^2))$, where ($a, b, c, d\in N$) and $c^2+cd+d^2=p$ is prime number; it is obvious that $(c, d) = 1$.

As it is clear from (115) and (116) is clear are possible the following cases. Let's consider each of them:

$$\text{I. } \begin{cases}ac-bd\equiv 0(\mathrm{mod}\,p),\\ ad-bc\equiv 0(\mathrm{mod}\,p),\\ c^2+cd+d^2=p,\\ a^2+ab+b^2\equiv 0(\mathrm{mod}\,p).\end{cases} \Rightarrow \begin{cases}(c-d)(a+b)\equiv 0(\mathrm{mod}(c^2+cd+d^2)),\\ a^2+ab+b^2\equiv 0(\mathrm{mod}(c^2+cd+d^2)).\\ (a,b)=1,\ (c,d)=1.\end{cases} \Rightarrow$$

$$\Rightarrow \begin{cases}a+b\equiv 0(\mathrm{mod}\,p),\\ a^2+ab+b^2\equiv 0(\mathrm{mod}\,p),\\ (a,b)=1\end{cases} \Rightarrow \begin{cases}a^2+2ab+b^2\equiv 0(\mathrm{mod}\,p),\\ a^2+ab+b^2\equiv 0(\mathrm{mod}\,p),\\ (a,b)=1\end{cases} \Rightarrow \begin{cases}ab\equiv 0(\mathrm{mod}\,p),\\ a^2+ab+b^2\equiv 0(\mathrm{mod}\,p),\Rightarrow\\ (a,b)=1\end{cases}$$

$$\Rightarrow \begin{cases}a\equiv 0(\mathrm{mod}\,p),\\ b\equiv 0(\mathrm{mod}\,p),\Rightarrow p=1.\ \text{this is impossible.}\\ (a,b)=1.\end{cases}$$

$$\text{II. } \begin{cases}ac-bd\equiv 0(\mathrm{mod}\,p),\\ ad+bc+bd\equiv 0(\mathrm{mod}\,p).\end{cases}$$ From (113) it is obvious that presentation is possible;

III. $\begin{cases} ac+bd+bc \equiv 0 (\bmod p), \\ ad-bc \equiv 0 (\bmod p). \end{cases}$ From (114) it is obvious that presentation is possible;

IV. $\begin{cases} ac+bd+bc \equiv 0 (\bmod p), \\ ad+bc+bd \equiv 0 (\bmod p), \\ a^2+ab+b^2 \equiv 0 (\bmod p), \\ (a,b)=1, (c,d)=1. \end{cases} \Rightarrow \begin{cases} a(c-d) \equiv 0 (\bmod (c^2+cd+d^2)), \\ a^2+ab+b^2 \equiv 0 (\bmod (c^2+cd+d^2)), \Rightarrow \\ (a,b)=1,\ (c,d)=1. \end{cases}$

$\begin{cases} a \equiv 0 (\bmod p), \\ a^2+ab+b^2 \equiv 0 (\bmod p), \Rightarrow \\ (a,b)=1. \end{cases} \begin{cases} a \equiv 0 (\bmod p), \\ b \equiv 0 (\bmod p), \Rightarrow p=1. \\ (a,b)=1. \end{cases}$ This is impossible.

Proceeding from all this, by dividing (113) and (114) on $(c^2+cd+d^2)^2$ we will obtain the provable.

I.e. if $a^2+ab+b^2 \equiv 0(\bmod(c^2+cd+d^2))$, where ($a,b,c$,d$\in N$), ($a,b$)=1) and $c^2+cd+d^2=p$ are prime numbers, then the dividend also will be as (i.e. it would be represented as $A^2+AB+B^2$, where $A$, $B \in N$. Q.E.D.

this form, the dividend has a divisor that cannot be represented in the form (118).

**Lemma 10**. If $a^2+ab+b^2$ (117) type number, where $a$, $b \in N$ and ($a$, $b$) = 1, is divisible by a prime numbers that would not be represented in such way, the dividend has a divisor that is presented as (118).

**Proof.**

Let's assume that

$$a^2+ab+b^2 = x \cdot P_1, P_2, \dots, P_n. \tag{119}$$

If $P_1, P_2, \dots, P_n$ all prime divisor will be as (118), then by successive division of (119) on $P_1$, $P_2$ and so on $P_n$ by virtue of the previous theorem, each divisor including $x$ must be as (118). And this contradicts the condition of Lemma 3. Therefore, any of $P_i$ $i=1,\dots,n$ does not have such form. Q.E.D.

Now let's prove

Theorem 5. If $a$, $b \in N$ and ($a$, $b$) = 1, then each divisors of $(a^2+ab+b^2)$ will be as (120).

Proof:

Let's say that ($a$, $b$) = 1, $a$, $b \in N$ and divisor $x$ of $(a^2+ab+b^2)$ is not as (120). In addition let's say that $a^2+ab+b^2$ is east between such numbers. It is obvious that $a$ and $b$ would be presented as following:

$$a = mx + (-1)^k c\,; \;\; b = nx + (-1)^\ell d\,, \tag{121}$$

where $m, n \in Z_0$; $c, d, k, \ell \in N$; $0 < x < \frac{x}{2}$, $0 < d < \frac{x}{2}$ .

then

$$a^2 + ab + b^2 = (mx + (-1)^k c)^2 + (mx + (-1)^k c)(nx + (-1)^\ell d) + (nx + (-1)^\ell d)^2 =$$

$$= (m^2 + mn + n^2)x^2 + (2(-1)^k mc + (-1)^\ell md + (-1)^k nc + 2(-1)^\ell nd)x +$$

$$+(c^2 + (-1)^{k+1} cd + d^2)\,. \tag{122}$$

Let's mention that if $k$ and $\ell$ simultaneously are odd, then $mn \neq 0$. (123)

Let's show that by conditions of (121) and (122) in the right side of (40) the sum of first two summands is positive. Indeed

$$(m^2 + mn + n^2)x^2 + (2(-1)^k mc + (-1)^\ell md + (-1)^k nc + 2(-1)^\ell md)x = (m^2 + mn + n^2)x^2 +$$

$$+(-2mc - md - nc - 2nd)x > (m^2 + mn + n^2)x^2 + \left(-2m\frac{x}{2} - m\frac{x}{2} - n\frac{x}{2} - 2n\frac{x}{2}\right)x =$$

$$= (2m^2 + 2mn + 2n^2 - 3m - 3n)\frac{x^2}{2} = (m(2m + n - 3) + n(2n + m - 3))\frac{x^2}{2} > 0\,, \tag{124}$$

As from $m \geq 1$; $n \geq 1$ and $x \in N$, i.e.

$$c^2 + (-1)^{k+\ell} cd + d^2 < a^2 + ab + b^2 \text{ and } c^2 + (-1)^{k+\ell} cd + d^2 \equiv 0 (\mathrm{mod}\, x)\,. \tag{125}$$

Let's mention that

$$c^2 - cd + d^2 = (c - d)^2 + d(c - d) + d^2, \tag{126}$$

$$c^2 - cd + d^2 = (d - c)^2 + d(d - c) + c^2, \tag{127}$$

$$\begin{cases} 0 < c < \frac{x}{2}, \\ 0 < d < \frac{x}{2}. \end{cases} \Rightarrow |c - d| < \frac{x}{2}. \tag{128}$$

Let's say that $k$ and $\ell$ are equally even or odd, then from (125) we have

$$c^2 + cd + d^2 = yx\,. \tag{129}$$

If $c$ and $d$ has more than 1 common divisors, then it is not dividing on $x$, as then it will divide also $a$ and $b$, but $(a, b) = 1$. By dividing of (129) on square of common divisor we will obtain

$$e^2 + ef + f^2 = zx\,, \tag{130}$$

In addition

$$e^2 + ef + f^2 < c^2 + cd + d^2 < \frac{3x^2}{4} < x^2 < a^2 + ab + b^2\,.$$

I.e. we obtain $e^2+ef+f^2<a^2+ab+b^2$ and it also has divisor that is not as (120). This would lead us to the infinite assumption.

I.e. if ($a$, $b$) = 1, then (120) has only and only (120) type divisors.

If $k$ and $\ell$ are various odd and even numbers then we have:

$c^2+cd+d^2=yx$ and if $c > d$, then by virtue of (126)

$$(c-d)^2+d(c-d)+d^2=yx, \tag{131}$$

And if $c > d$, then from (127)

$$(d-c)^2+c(d-c)+c^2=xy. \tag{132}$$

By similar reasoning, as for (129), we will obtain the proof or if ($a$, $b$) = 1, (120) has only $(a^2+ab+b^2)$ type divisors.

Application of this theorem is possible in lot of tasks. Let's state one example. Let's solve the equation $x^2+xy+y^2=n$. Or prior solving this equation, we can first determine whether this equation is solvable. In the canonic expansion of $n$ from prime numbers with odd powers, if any of them does not have the first form, then this equation has no solution.

Here we would like to briefly considered our studies in scope of Fibonacci numbers.

A lot of works has been written on Fibbonacci numbers and its properties should probably be almost completely studied, but as Morley's theorem, the magnificent property that any triangle at the intersection of trisections of any angle (both internal and external angles) will always turn out to be an equilateral triangle, and from the discovery of that less than 90 years has passed, it makes one wonder day it is possible to discover Fibonacci numbers with currently unknown amazing properties.

For our part, we want to draw attention to several facts.

It is known that Bines, when deriving his formula

It is known that Bines at deriving his formula

$$u_n=\frac{\left(\frac{1+\sqrt{5}}{2}\right)^n-\left(\frac{1-\sqrt{5}}{2}\right)^n}{\sqrt{5}}, \tag{133}$$

Is applied the equation

$$x^2=x+1, \tag{134}$$

The question is raised why exactly the equation (134)? Formula (134) itself raises the question of the expediency of considering the following equations. These equations are:

$$x^3=x^2+x+1; \tag{135}$$

$$x^4=x^3+x^2+x+1; \tag{136}$$

$$x^5=x^4+x^3+x^2+x+1, \text{ (137)} \qquad \ldots; x^n=x^{n-1}+x^{n-2}+\ldots+x^2+x+1. \tag{138}$$

The diagrams of sides of equation (135) are easily understandable for us, and even for (136) after much torment it will become more or less clear. And we would lose all interest in the subsequent equations.

The situation is different if we proceed as follows:

$$\begin{cases} n \in N; \\ x^n = x^{n-1} + x^{n-2} + \ldots + x^2 + x + 1; \\ x \ge 1. \end{cases} \begin{cases} n \in N; \\ x^n(x-1) = (x-1)(x^{n-1} + n^{n-2} + \cdots + x^2 + x + 1); \Leftrightarrow \\ x \ge 1 \end{cases}$$

$$\Leftrightarrow \begin{cases} n \in N; \\ x^n(x-1) = x^n - 1; \\ x \ge 1. \end{cases} \Leftrightarrow \begin{cases} n \in N; \\ x^n = \dfrac{1}{2-x}; \\ x \ge 1, \quad 2 - x > 0. \end{cases} \Leftrightarrow \begin{cases} n \in N; \\ x^n = \dfrac{1}{2-x}; \\ x \in [1; 2). \end{cases}$$

Now it becomes clear what we are dealing with.

Let's consider $x^n = \dfrac{1}{2-x}$ $\quad n = 1 \Leftrightarrow x = 1$

$$n = 2 \Rightarrow \left[ \begin{array}{l} x = \dfrac{1-\sqrt{5}}{2} \\ x = 1 \\ x = \dfrac{1+\sqrt{5}}{2} \approx 1,6180 \end{array} \right. \qquad n = 3 \Rightarrow \left[ \begin{array}{l} x = 1 \\ x = 1,83929 \end{array} \right. \qquad n = 4 \Rightarrow \left[ \begin{array}{l} x = -0,7748 \\ x = 1 \\ x = 1,9276 \end{array} \right.$$

$$n = 5 \Rightarrow \left[ \begin{array}{l} x = 1 \\ x = 1,96595 \end{array} \right. \qquad n = 6 \Rightarrow \left[ \begin{array}{l} x = -0,8403 \\ x = 1 \\ x = 1,98358 \end{array} \right. \qquad n = 7 \Rightarrow \left[ \begin{array}{l} x = 1 \\ x = 1,99196 \end{array} \right.$$

$$n = 8 \Rightarrow \left[ \begin{array}{l} x = -0,87629 \\ x = 1 \\ x = 1,99603 \end{array} \right. \qquad n = 9 \Rightarrow \left[ \begin{array}{l} x = 1 \\ x = 1,99803 \end{array} \right.$$

Here we have obtained interesting characteristic numbers. The number that characterizes the golden ratio $q = 1,6180$ and it is probably necessary to pay attention to the obtained values: 1; 1.83929; 1.9276; 1.96595; 1.98358; 1.99603; 1.99196; 1.99803 and so on.

In a few words, let's touch on some construction problems.

Let's say that is given any *MN* segment, $|MN| = a \ne 1$, then in order to construct segments with lengths: $a^2$; $a^3$; $a^4; \cdots; a^n$; , is necessary to do the following.

1) Let's on arbitrary straight line *p* take the points *D* and *B* in such manner that: $|DB| = 1$, $B\hat{D}C = 90°$ and $|CD| = a$. Then if from *C* point we draw a perpendicular to *CB* segment up to intersection of circumference *p* , we will obtain that

$$|CD|^2 = |AD||BD| \Rightarrow a^2 = |AD|\cdot 1 \Rightarrow |AD| = a^2 .$$

2) Then similarly let's construct $\Delta A_1C_1B_1$, where $|D_1B_1| = a$; $|C_1D_1| = a^2$;

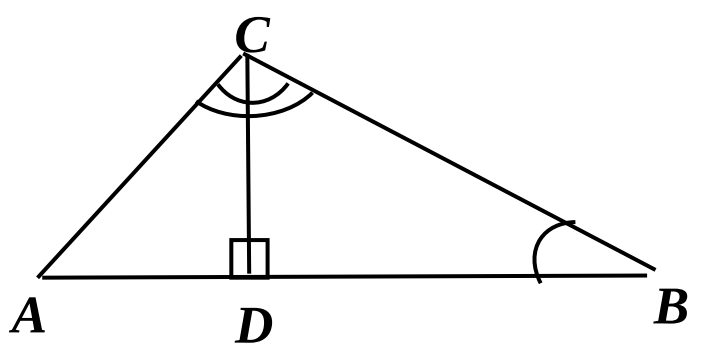

$A_1\hat{C}_1B_1 = C_1\hat{D}_1B_1 = 90°$, then similarly $|C_1D_1|^2 = |A_1D_1|\cdot|D_1B_1| \Rightarrow$ $\Rightarrow (a^2)^2 = |A_1D_1|\cdot a \Rightarrow |A_1D_1| = a^3$ and so on.

Let's construct: $\Delta A_{n-1}C_{n-1}B_{n-1}$, where

$$|D_{n-1}B_{n-1}| = a^{n-2};\ |C_{n-1}D_{n-1}| = a^{n-1};\ A_{n-1}\hat{C}_{n-1}B_{n-1} = C_{n-1}\hat{D}_{n-1}B_{n-1} = 90° .$$

Here $|C_{n-1}D_{n-1}|^2 = |A_{n-1}D_{n-1}|\cdot|D_{n-1}B_{n-1}| \Rightarrow (a^{n-1})^2 = |A_{n-1}D_{n-1}|\cdot a^{n-2} \Rightarrow |A_{n-1}D_{n-1}| = a^n$.

I.e. if is given any $[MN]$, such that $|MN| = a \neq 1$, then is easy to construct the segment with length of $a^n$. It is known that also is easy to construct the segment with length of (2 – $a$).

Let's consider the following sequences:

1) $u_1 = u_2 = u_3 = 1$, $\quad u_n = u_{n-1} + u_{n-2} + u_{n-3}$;

2) $u_1 = u_2 = u_3 = u_4 = 1$, $\quad u_n = u_{n-1} + u_{n-2} + u_{n-3} + u_{n-4}$

3) $u_1 = u_2 = \ldots u_k = 1$, $\quad u_n = u_{n-1} + \ldots + u_{n-k}$.

1) Pseudo-Fibbonacci numbers are sequences of order 3; 2) Pseudo-Fibbonacci numbers are sequences of order 4; 3) sequences of order $k$. Let's call these sequences as pseudo-Fibbonacci numbers: Let's call these sequences as pseudo-Fibbonacci numbers: 1) a sequence of 3 orders is pseudo-Fibbonacci numbers; 2) pseudo-Fibbonacci numbers is a sequence of 4 orders; 3) order sequence.

Now consider the first sequence

Let's now consider the Let's call these sequences pseudo-Fibbonacci numbers: 1) a sequence of 3 orders is pseudo-Fibbonacci numbers; 2) pseudo-Fibbonacci numbers is a sequence of 4 orders; 3) order sequence.

Now consider the first sequence

$u_1 = u_2 = u_3 = 1$, $u_4 = 3$, $u_5 = 5$, $u_6 = 9$, $u_7 = 17$, $u_8 = 31$, $u_9 = 57$,

$u_{10} = 105$, $u_{11} = 193$. It is easy to show that Let's prove

$$\boxed{u_n = 2u_{n-1} - u_{n-4}} \tag{139}$$

$S_1 = 1$, $S_2 = 2$, $S_3 = 3$, $S_4 = 6$, $S_5 = 11$, $S_6 = 20$, $S_7 = 37$, $S_8 = 68$.

$$\boxed{S_n = 2S_{n-1} - S_{n-4}} \tag{140}$$

It is easy to notice that

$$u_5 = 2u_4 - 1 = 2u_4 - u_1$$

$$u_6 = 2u_5 - 1 = 2u_5 - u_2$$
$$u_7 = 2u_6 - 1 = 2u_6 - u_3$$
$$u_8 = 2u_7 - 3 = 2u_7 - u_4$$
$$u_9 = 2u_8 - 5 = 2u_8 - u_5$$
$$u_{10} = 2u_9 - 9 = 2u_9 - u_6$$
$$\ldots\ldots\ldots\ldots\ldots\ldots\ldots\ldots\ldots\ldots\ldots$$
$$u_n = 2u_{n-1} - u_{n-4} \quad \text{by summation.}$$

$$u_5 + u_6 + \ldots + u_n = 2(u_4 + u_5 + \ldots + u_{n-1}) - (u_1 + u_2 + u_3 + \ldots + u_{n-4}) \Rightarrow$$
$$\Rightarrow 6 + u_5 + u_6 + \ldots + u_n = 2(3 + u_4 + u_5 + \ldots + u_{n-1}) - (u_1 + u_2 + \ldots + u_{n-4}) \Rightarrow$$
$$\Rightarrow \ S_n = 2S_{n-1} - S_{n-4} \ \text{ Q.E.D.}$$

It would be easily shown that if $S_1 + S_2 + \ldots + S_n \equiv S'_n$, then similarly we will obtain $S'_n = 2S'_{n-1} - S'_{n-4}$ for each $n \in N \setminus \{1; 2; 3; 4\}$ and so infinitely.

Pseudo-Fibbonacci numbers were determined recursively, individually, according to their numbers. It turns out that we would directly find any term of this sequence.

Let's study different sequences $u_1, u_2, \ldots, u_n, \ldots$ that satisfy the dependency

$$u_n = 2u_{n-1} - u_{n-4} \ . \qquad (141)$$

All such sequences are called by us as solutions of equation (141).

Let's accordingly designate them as $V$, $V'$, $V''$ and $V'''$ sequences

$$v_1,\ v_2,\ v_3, \cdots$$
$$v'_1,\ v'_2,\ v'_3, \cdots$$
$$v''_1,\ v''_2,\ v''_3, \cdots$$
$$v'''_1,\ v'''_2,\ v'''_3, \cdots$$

Firstly let's prove two lemmas.

Lemma 1. If $V$ is the solution of (141), and $c \in R$, then $cV$ (or the sequence $cv_1, cv_2, cv_3, \ldots$) also is the solution of (141).

Let's multiply the dependency $v_n = 2v_{n-1} - v_{n-4}$ term wise $c$, we will obtain

$$cv_n = 2cv_{n-1} - cv_{n-4} \ \text{ Q.E.D.}$$

Lemma 2. If $V'$, $V''$ and $V'''$ are the solutions of (141), then their sum $V' + V'' + V'''$ (or the sequence $v'_1 + v''_1 + v'''_1$, $v'_2 + v''_2 + v'''_2$, $v'_3 + v''_3 + v'''_3, \ldots$) also represents the solution of (141).

From the conditions of lemma we have

$$v'_n = 2v'_{n-1} - v'_{n-4}$$

and

$$v''_n = 2v''_{n-1} - v''_{n-4}\,,\;\; v'''_n = 2v'''_{n-1} - v'''_{n-4}\;.$$

Due term wise summation of this equalities we will obtain

$$v'_n + v''_n + v'''_n = 2(v'_{n-1} + v''_{n-1} + v'''_{n-1}) - (v'_{n-4} + v''_{n-4} + v'''_{n-4})\,.$$

Due this Lemma is proven.

Let's say that now $V'$, $V''$, $V'''$ are three non-proportional solutions of (141) (or are two such solutions of (141), that for arbitrary $c$ would be found such number *n*, for that $\frac{v'_n}{v''_n} \neq c$). Let's show that each sequence *V* that represents the solution of (141), is possible to represent as

$$c_1V' + c_2V'' + c_3V''',\qquad (142)$$

where $c_1$ and $c_2$ – are certain constants. Therefore, it is customary to say that (142) represents the general solution of (141).

Now let's take certain solution *V* of (141). As already was mentioned, a sequence is considered as given if its first three members $v_1$, $v_2$ and $v_3$ are given.

Let's find such $c_1$, $c_2$ and $c_3$, for that

$$\begin{cases} c_1v'_1 + c_2v''_1 + c_3v'''_1 = v_1 \\ c_1v'_2 + c_2v''_2 + c_3v'''_2 = v_2 \\ c_1v'_3 + c_2v''_3 + c_3v'''_3 = v_3 \end{cases}\qquad (143)$$

From the condition of (141) is solvable the system with respect of (143) $c_1$, $c_2$ and $c_3$

$$c_1 = \frac{\Delta c_1}{\Delta}\quad \Delta c_1 = \begin{vmatrix} v_1 & v''_1 & v'''_1 \\ v_2 & v''_2 & v'''_2 \\ v_3 & v''_3 & v'''_3 \end{vmatrix}\qquad c_2 = \frac{\Delta c_2}{\Delta}\quad \Delta c_2 = \begin{vmatrix} v'_1 & v_1 & v'''_3 \\ v'_2 & v_2 & v'''_2 \\ v'_3 & v_3 & v'''_3 \end{vmatrix}$$

$$c_3 = \frac{\Delta c_3}{\Delta}\quad \Delta c_3 = \begin{vmatrix} v'_1 & v''_1 & v_1 \\ v'_2 & v''_2 & v_2 \\ v'_3 & v''_3 & v_3 \end{vmatrix}\qquad \Delta = \begin{vmatrix} v'_1 & v''_1 & v'''_1 \\ v'_2 & v''_2 & v'''_2 \\ v'_3 & v''_3 & v'''_3 \end{vmatrix}.$$

The condition (141) mean that $\Delta \neq 0$. If we introduce the obtains $c_1$, $c_2$ and $c_3$ in (142), we obtain the desired sequence *V*.

Therefore, to obtain all solutions of equation (141), is sufficient to find its three non-proportional solutions.

Let's start the search among geometric progressions. According to Lemma 2, it is sufficient to limit by the consideration of such progressions, the first term of that is 1. Thus, we take the progression

$$1,\; q,\; q^2,\; q^3, \ldots$$

In order for it to be a solution (141) for this, it is necessary for each $\begin{cases} n > 4 \\ n \in N \end{cases}$ to fulfill

$$q^n = 2q^{n-1} - q^{n-4},$$

Or as $q \neq 1$, we obtain

$$q^3 - q^2 - q - 1 = 0. \qquad (144)$$

Let's say that solutions of (144) are $q_1$, $q_2$, $q_3$, then

$$q_1 = \frac{\sqrt[3]{19+3\sqrt{33}} + \sqrt[3]{19-3\sqrt{33}} + 1}{3}$$

$$q_2 = -\frac{\sqrt[3]{19+3\sqrt{33}} + \sqrt[3]{19-3\sqrt{33}} - 2}{6} + i\sqrt{3}\frac{(\sqrt[3]{19+3\sqrt{33}} - \sqrt[3]{19-3\sqrt{33}})}{6}$$

$$q_3 = -\frac{\sqrt[3]{19+3\sqrt{33}} + \sqrt[3]{19-3\sqrt{33}} - 2}{6} - i\sqrt{3}\frac{(\sqrt[3]{19+3\sqrt{33}} - \sqrt[3]{19-3\sqrt{33}})}{6}.$$

Let's mention that for each root would be fulfilled

$$1 + q_1 + q_1^2 = q_1^3, \quad 1 + q_2 + q_2^2 = q_2^3, \quad 1 + q_3 + q_3^2 = q_3^3 \quad \text{and} \quad q_1 q_2 q_3 = -1.$$

We obtain three geometric progressions, that represents the solutions of (141). Thus the progression

$$c_1 + c_2 + c_3, \ c_1 q_1 + c_2 q_2 + c_2 q_2, \ c_1 q_1^2 + c_1 q_2^2 + c_3 q_3^2, \qquad (145)$$

Also represents the solution of (141). As the found progressions have different values, they are therefore disproportionate. Formula (142) for various $c_1$, $c_2$ and $c_3$ will give all solutions of (141).

In particular, for any $c_1$, $c_2$ and $c_3$ formula (145) gives a pseudo-Fibbonacci series. To do this, we rewrite (145) in this way .

$$\begin{cases} c_1 + c_2 + c_3 = 1 \\ c_1 q_1 + c_2 q_2 + c_3 q_3 = 1 \\ c_1 q_1^2 + c_2 q_2^2 + c_3 q_3^2 = 1 \end{cases}$$

$$c_1 = \frac{\Delta c_1}{\Delta}, \quad c_2 = \frac{\Delta c_2}{\Delta}, \quad c_3 = \frac{\Delta c_3}{\Delta}.$$

$$\Delta c_1 = \begin{vmatrix} 1 & 1 & 1 \\ 1 & q_2 & q_3 \\ 1 & q_2^2 & q_3^2 \end{vmatrix}, \quad \Delta c_2 = \begin{vmatrix} 1 & 1 & 1 \\ q_1 & 1 & q_3 \\ q_1^2 & 1 & q_3^2 \end{vmatrix}, \quad \Delta c_3 = \begin{vmatrix} 1 & 1 & 1 \\ q_1 & q_2 & 1 \\ q_1^2 & q_2^2 & 1 \end{vmatrix}, \quad \Delta = \begin{vmatrix} 1 & 1 & 1 \\ q_1 & q_2 & q_3 \\ q_1^2 & q_2^2 & q_3^2 \end{vmatrix}.$$

$$c_3 = \frac{(1-q_1)(1-q_2)}{(q_3-q_1)(q_3-q_2)}, \quad c_2 = \frac{(1-q_1)(q_3-1)}{(q_2-q_1)(q_3-q_2)}, \quad c_1 = \frac{(1-q_3)(1-q_2)}{(q_3-q_1)(q_2-q_1)}.$$

$$u_n = c_1 q_1^{n-1} + c_2 q_2^{n-1} + c_3 q_3^{n-1}.$$

Similarly we behave in the next case when

$u_1 = u_2 = u_3 = u_4 = 1$ and so on.

**Conclusions**

We have shown that the entire class of tasks would be easily solved using one simple lemma; in addition, is stated the task to study $\frac{p^q - 1}{p - 1}$ type expressions, where $p$ and $q$ are prime numbers, with properties; some properties of the $(a^2 + ab + b^2)$ type expression were studied, and so on.

The following issues are considered in the work:

– Obviously, for prime numbers $p$ and $q$, it is of great interest to determine the quantity of those prime divisors of the $A = \frac{p^q - 1}{p - 1}$ number that are less than $p$. To do this, we have considered:

Theorem 1. Let's say that $p$ and $q$ are odd prime numbers and $p = 2q + 1$. Then from the various prime divisors of number $A = \frac{p^q - 1}{p - 1}$, taken separately, only one is less than $p$. $A$ has at least two different prime divisors.

Theorem 2. Let's say that $p$ and $q$ are odd prime numbers and $p < 2q + 1$. Then all prime divisors of the number $A = \frac{p^q - 1}{p - 1}$ are greater than $p$;

Theorem 3. Let's say that $q$ is an odd prime number and $p \in N \setminus \{1\}$, $p \in ]1; q] \bigcup [q + 2; 2q]$, then each of the various prime divisors of the number $A = \frac{p^q - 1}{p - 1}$, taken separately, is greater than $p$;

Theorem 4. Let's say that $q$ is a prime odd number and $p \in \{q + 1; 2q + 1\}$, then from different prime divisors of number $A = \frac{p^q - 1}{p - 1}$, taken separately, only one of them is less than $p$. $A$ has at least two different prime divisors.

Task 1. Let's solve the equation $2^x = \frac{y^z - 1}{y - 1}$ in natural $x$, $y$, $z$ numbers. At the same time $y$ must be a prime number.

Task 2. Let's solve the equation $3^x = \frac{y^z - 1}{y - 1}$ in natural $x$, $y$, $z$ numbers. At the same time $y$ must be a prime number.

Task 3. Let's solve the equation $p^x = \frac{y^z - 1}{y - 1}$, where $p \in \{5; 7; 11; 13; ...\}$ prime number $x$, $y \in N$ and $y$ – is prime number.

I – is stated the Lemma, by that easily would be solved the class of tasks:

Lemma 1. Let's say that $a$, $b$, $n \in N$ and $(a, b) = 1$. Let's prove that if $a^n \equiv 0 \pmod{|a-b|}$, or $b^n \equiv 0 \pmod{|a-b|}$, then $|a-b| = 1$.

Let's solve the equations (I – X) in natural $x$, $y$ numbers:

I. $\left(\frac{x+y}{2}\right)^z = x^z - y^z$;

II. $(x+y)^z = (2x)^z + y^z$;

III. $(x+y)^z = (3x)^z + y^z$;

IV. $(y-x)^{x+y} = x^y$, $(y > x)$;

V. $(y-x)^{x+y} = y^x$, $(y > x)$;

VI. $(x+y)^{x-y} = x^y$;

VII. $(x+y)^{x-y} = y^x$;

VIII. $(x+y)^y = (x-y)^x$, $(x > y)$;

IX. $(x-y)^{x+y} = x^{x-y}$;

X. $(x+y)^{x-y} = (x-y)^x$, $(y > x)$.

theorem 5. If a, $b \in N$, $(a, b) = 1$, then each divisor of $(a^2 + ab + b^2)$ will be as similar.

## Refferences

1. Aghdgomelashvili Zurab. Diophantine geometric figures (Diophantine, Bidiophantine, Pseudodiophantine and Pseudodiophantine planar geometric (Figures) LAP LAMBERT Academic Publishing. 2020 (978-720-2-52393-6).
2. Aghdgomelashvili Zurab. On a Fundamental Task of Diophantine Geometric Figures. Cornell University. Arxiv > org>arxiv. 2003. 02652.
3. Aghdgomelashvili Zurab. On a Fundamental Task of Bidiophantine Geometric Figures. Cornell University. Ardxiv > org>arxiv. 2003. 10 846.
4. Zurab Agdgomelashvili. About on fundamental task on diophantine geometric Figures. Georgian Technical University. N 4(514), 2019, 141-173.
5. Zurab Agdgomelashvili. Some interesting tasks from the classical number theory. Georgian Technical University. N 4(518), 2020, 150-188.
6. Matiyasevich Y, Hilbert's Tenth Problem. MLT Press Cambridge, Massachusetts, 1993.
7. Yaglom A.M. and Yaglom I.M. Non elementary tasks in the elementary statement. State Publishing of Technical-theoretical literature. Moscow, 1954.